\documentclass{amsart}
\usepackage{amsmath,color}
\usepackage{amssymb}
\usepackage{amsfonts}
\usepackage{amsthm}
\usepackage{amscd}

\newtheorem{thm}{Theorem}[section]

\newtheorem{defn}[thm]{Definition}
\renewcommand{\proof}[1][Proof]{\noindent\textsc{#1}. }

\begin{document}

\title{ Musical modes, their associated chords and their musicality}

\author{Mihail Cocos \& Kent Kidman}

\begin{abstract}

In this paper we present a mathematical way of defining musical modes and we define the musicality of a mode as a product of three different factors.  We conclude by classifying the modes which are most musical according to our definition.

\end{abstract}

\maketitle

\section{Musical scales and modes} 
 First we must define some concepts.  In traditional western style of music there are $12$ different notes in an octave.  Two notes represent the same tone if they are at a distance of $12n$ semitones (a multiple of octaves apart).  The $13^{th}$ note restarts the cycle and it is the same tone as our first note (simply an octave higher).  Loosely speaking a scale is simply an ordered sequence of different tones within the same octave.  For instance $CDEFGAB$ is the well known C major scale, $ ABCDEFG$ is the A minor scale, and $EGAA^{\#}BD$ is the hexatonic blues scale in the key of E.  The first note of the scale is said to be the key of the scale.\\

\begin{defn}

  For further considerations about scales we will associate every tone with a numerical value, namely its position in the chromatic A scale \[AA^{\#}BCC^{\#}CDD^{\#}EFF^{\#}GG^{\#}\]
The lexicographic order of tones will be defined as the order induced by their numerical values defined as above.\\

\end{defn}

\noindent The precise mathematical definition of a scale is:

\begin{defn}
A $k$ scale is a circular permutation of $k$ different tones in lexicographic order.

\end{defn}
\noindent To exemplify the definition let's take $7$ different tones in increasing order, namely $ ABCDEFG.$  This is simply the A minor scale.  The next circular permutation of the same tones is $BCDEFGA$ and is called B Mixolydian.  The following is $CDEFGAB$ and is called C major, and thus circularly permuting the notes of the A minor scale we obtain $7$ different types of scales used in western music.  It is well known that these scales have the same number of chords (harmonic triads).
A song composed in the A minor key simply uses only the tones of the scale A minor.  For instance a song written in  the key of A minor will use the same tones as a song in the C major key.  In order to understand what a musical mode is, let's examine the following two scales.  Take the first one to be $CDEFGAB$ (C major) and the second one $ABC^{\#}DEF^{\#}G^{\#}$(A major).  Although they start on a different note when played they sound quite similar.  The reason is that the number of semitones between two consecutive notes in both of them follow the same pattern, that is classically known as w,w,h,w,w,w,h (where w=whole= 2 semitones and h=half= 1 semitone).  Two scales that have a similar interval pattern are said to be of the same mode (thus the name major for both of the scales considered above).  Now assume that a scale starts on the note D and has the same interval pattern (major).  It's easy to see that the notes of the D major scale will be $ DEF^{\#}GABC^{\#}$, and thus we obtain another scale in the major mode.  Any song that is played in a certain major scale can be played in any other major scale and the two scales will have an equal number of chords.  In what follows we will make more clear (mathematical) the notions of a mode and chord and we will classify modes by their musicality (a certain product of factors that we believe make the musicality of a scale).
 \begin{defn}

  A $k$ ( $k>2$) mode is a $k-$tuple of positive numbers \[ (a_1,a_2, \cdots , a_k) \] such that  \[\sum_{i=1}^{k} a_i  = 12.\]

\end{defn}

\noindent The first $(k-1)$ numbers in the $k$ tuple represents the absolute value of the difference of the numerical values between two consecutive notes in a given scale, while the last represents the left over number of semitones between the last note of the scale and the represent-ant of the first one in the next octave.
The reason behind this definition is that musical literature has exotic scales with more or less than $7$ different notes and for which the interval between two consecutive notes is more than a whole step.  As an example consider the blues scale  based on E mentioned before, namely $ EGAA^{\#}BD$. Its interval pattern (i.e mode) is $(3,2,1,1,3,2)$.  Another exotic well known example is the pentatonic  In Sen scale based in A, that is $ AA^{\#}DEG$.  Its mode is $(1,4,2,3,2).$  We have to note the fact that this is a purely mathematical definition and not every single scale obtained by  choosing a random mode will sound musical when played. For instance if we choose the mode of a scale to be $ (1,1,1,1,1,1,1,1,1,1,1,1) $ and its key to be A then we obtain the chromatic scale based in A that consists of $ 13 $ half steps.  This is NOT a very musical scale.

\begin{thm}

The number of $ k $ modes is  $  \left (  \begin{array}{c} 11 \\
 k-1 \end{array}  \right ).  $

\end{thm}

\proof The number of $k$ modes is obviously equal to the number of scales with $k$ tones based in A.  It follows that counting  all possible scales with $k$ different tones will yield the number of $k$ modes. To form a scale based in A we need to choose the leftover $k-1$ tones (first one being A is already chosen) out of the remaining $11$ tones. Thus the number of $k$ tones is  $  \left (  \begin{array}{c} 11 \\
 k-1 \end{array}  \right )  .$                  $\Box$

\bigskip

In order to give a precise definition of the musicality of a mode we will need to define the notion of a harmonic chord. Here we have to recall that choosing a key and a mode we obtain a scale. We will give our definition of a harmonic triad (chord) starting from a scale. It will be clear from our definition that a change of the key will only change the name of the chords and not their number.

\begin{defn}  A harmonic interval is an interval of $3, 4, 5, 7, 8,$ or $9$ semitones. 

\end{defn}

\noindent {\bf Note:}  Western music has evolved to where these intervals are considered generally more pleasant to the ear than other intervals between notes.  These intervals form the building blocks for major and minor chords. Jazz music and often times even classical also consider an interval of $6$ semitones to be an acceptable interval for the formation of a chord. 

\medskip

\noindent We will adopt and adapt the next two definitions and theorem from \cite{fowers}.	
 	
\begin{defn}  A harmonic subset of a scale $S$ is a subset of notes in $S$ such that the interval between any two of its notes is a harmonic interval.
\end{defn}

\begin{thm}
 The maximum number of notes in a harmonic subset is three.
 \end{thm}

\proof Let $1\leq n_1<n_2< \cdots <n_s \leq 12,$ be the ordered numerical values of the notes in our harmonic subset and assume that $s>3.$ By the very definition of a harmonic subset $|n_i-n_j| \in \{3,4,5,7,8,9\}$ and looking at the first $4$ numerical values of the tones within our set $n_1<n_2<n_3<n_4$ we conclude that $n_4-n_1=9.$ It follows that $n_2-n_1=3, n_3-n_1=3,$ and $n_4-n_3=3,$ which in turn implies
that $n_3-n_1=6 \notin \{3,4,5,7,8,9\},$ which contradicts our hypothesis. $\Box$ \\

\noindent The above theorem allows us to give the following definition.

\begin{defn}

A $3$ tone harmonic subset of a scale is called a chord.

\end{defn}
\noindent {\bf Remark:}  According to the previous theorem (and here we're assuming that the harmonic intervals does not include the $6$ semitone intervals) it is clear that a chord consists of only $3$ different tones. If we include the $6$ semitone interval in our definition of harmonicity then the reader can easily verify that some chords can have a maximal number of $4$ different tones! Since these type of chords are simply two $3$ toned chords glued together we will not consider them  further in the study of musicality. Here is a list of the only types of harmonic chords $ n_1<n_2<n_3<n_4$:

\begin{itemize}
\item $n_2-n_1=3$ \& $n_3-n_2=4$ (major)
\item $n_2-n_1=4$ \& $n_3-n_2=3$ (minor)
\item $n_2-n_1=4$ \& $n_3-n_2=4$ (augmented)
\item $n_2-n_1=3$ \& $n_3-n_2=6$(minor sixth)
\item $n_2-n_1=3$ ,$n_3-n_2=3,$ \& $n_4-n_3=3$ (diminished)
\end{itemize}

\bigskip

The following statement is the motivation of the classification of modes in the next section.

\begin{thm}
If two scales belong to the same mode then they have the same number of chords and same number of one semi-tone intervals. Moreover if a mode is obtained as a circular permutation of another mode, the two modes will have the same number of chords and semi-tones.

\end{thm}

\noindent For a proof of this see \cite{fowers}. The previous Theorem will allow us, in the next section, to give a mathematical definition the definition of the "musicality" of mode.  This also leads to the following definition.
\bigskip

\begin{defn}

Two k modes are equivalent if one is a circular permutation of the other.

\end{defn}

\section{Classification of musical modes by their musicality}
In this section we present an attempt of defining the "musicality" of a mode. We will consider two modes to be equivalent if one is a circular permutation of the other one. 
\begin{defn}
The chord density of a $k$ mode, for $ k \ge 3$, is the ratio between the number of chords of a scale within that mode and $ \left (  \begin{array}{c} k \\
 3 \end{array}  \right )  $ and it is denoted $\mathcal{D}.$

\end{defn}

\noindent The motivation behind this definition is that most of the people will consider a scale to be more musical if it has more chords. In general a scale with more notes will likely have more chords, so in order to compare the chord richness of two scales with different number of notes we need to normalize the number of harmonic triads by the number of all possible combinations of three distinct notes  within the scale.

\begin{defn}
The tonal variety of a mode is $\mathcal{V}=\frac{k \sqrt{k}}{10},$ where $k$ the number of notes in the mode .  
\end{defn}
\noindent Obviously the more notes you have available to write a piece of music, the better the music will sound. So the "musicality" of a scale  should be direct proportional to the number of different notes. If we choose $v=\frac{k}{10}$ the most musical scales will be the ones with  only $3$ different tones, which is a fact that contradicts our intuition. If we choose $v=\frac{k^2}{10}$, then the most musical will be the ones with $11$ notes, again something that contradicts our intuition about musicality.

\begin{defn}
The anti-chromaticism of the mode is the number $\mathcal{A}$ of intervals larger than a semitone within the mode.

\end{defn}
\noindent The motivation of "anti-chromaticism" is that in general the most pleasantly sounding scales will have less semitone intervals.

\begin{defn}

If a mode has tonal variety $\mathcal{V}=\frac{k\sqrt{k}}{10}$, chord density $\mathcal{D}$, and anti-chromaticism $\mathcal{A}$, then its musicality is $ \mathcal{M}=\mathcal{D} \cdot \mathcal{V} \cdot \mathcal{A}.$
\end{defn}

We wish to find the $k$-modes which are most musical according to our definition, and a computer program was written to do this.  The program will search through all of the possible scales, count the number of chords in this scale, and keep track of those scales with the maximum number of chords.  The program will also check to only keep a representative scale of the set of all equivalent scales for a particular mode, equivalence being defined to be the same mode or a circular permutation of a mode (as described above).  In the end we will only have a few nonequivalent scales which are representatives of the modes with the maximum number of chords, and this will also give us the musicality of the mode as defined above.  To find all the combinations (scales), we used a modification of the algorithm found in \cite{Rosen}, p. 385.
\newpage
Here are the results of our algorithm:
\begin{table}[h!]
\caption{Musicality chart without the $6$ semitone interval}
\begin{tabular}{|l|l|l|l|l|}

\hline
Number of of tones&Scale representation&Number of chords&$\mathcal{A}$&Musicality\\
\hline
3 & 1 4 8   &  1 & 3 & 1.56\\
3 & 1 4 9   &  1 & 3 & 1.56\\
3 & 1 5 9   &  1 & 3 & 1.56\\
4 & 1 3 6 10   &  2 & 4 & 1.6\\
5 & 1 2 4 6 9   &  3 & 4 & 1.34\\
5 & 1 2 5 6 9   &  4 & 3 & 1.34\\
5 & 1 2 5 6 10   &  4 & 3 & 1.34\\
5 & 1 2 5 7 10   &  3 & 4 & 1.34\\
5 & 1 2 5 8 10   &  3 & 4 & 1.34\\
5 & 1 2 6 9 11   &  3 & 4 & 1.34\\
6 & 1 2 5 6 9 10   &  8 & 3 & 1.76\\
7 & 1 2 4 6 7 9 11   &  6 & 5 & 1.59\\
8 & 1 2 3 5 6 8 10 11   &  9 & 4 & 1.45\\
9 & 1 2 3 5 6 7 9 10 11   &  15 & 3 & 1.45\\
10 & 1 2 3 4 5 6 7 9 10 11   &  17 & 2 & .90\\
11 & 1 2 3 4 5 6 7 8 9 10 11   &  21 & 1 & .47\\
12 & 1 2 3 4 5 6 7 8 9 10 11 12   &  28 & 0 & 0.0\\

\hline
\end{tabular} 
\label{table1}
\end{table}

\bigskip

\noindent Let's us examine more closely the results of our algorithm displayed in the previous table. If we look at the (only!) $7$ notes scale in the table and we choose the key to be B then we get our familiar B  Mixolydian, which is a circular permutation of the C major scale. The scale with $6$ notes is a circular permutation of a variant of the blues scale (one that is more musical than the blues scale but harder to be played) and the 1,2,6,9,11 scale is only a variant of the In Sen scale. Thus we can see that our definition of musicality of a mode is consistent to the musicality as perceived by most of the people.
\newpage

\section{Musicality of modes with alternate harmonic definition}

It might be interesting to define a harmonic interval in a different way and examine mode musicality with a new definition.

\begin{defn}  An alternate harmonic interval is an interval of $3, 4, 5, 6, 7, 8,$ or $9$ semitones. 

\end{defn}

\noindent The only change we are making is adding the 6 semitone interval.  This is used in some jazz and classical music, as mentioned above.  By using this definition, we can get a new table of maximum musicality of nonequivalent modes.
\begin{table}[h!]
\caption{Musicality chart with the $6$ semitone interval}
\begin{tabular}{|l|l|l|l|l|}
\hline
Number of of tones&Scale representation&Number of chords&$\mathcal{A}$&Musicality\\
\hline
3 & 1 4 7   &  1 & 3 & 1.56\\
3 & 1 4 8   &  1 & 3 & 1.56\\
3 & 1 4 9   &  1 & 3 & 1.56\\
3 & 1 5 9   &  1 & 3 & 1.56\\
4 & 1 4 7 10   &  4 & 4 & 3.2\\
5 & 1 2 4 7 10   &  5 & 4 & 2.24\\
5 & 1 2 5 8 11   &  5 & 4 & 2.24\\
6 & 1 2 4 5 7 10   &  7 & 4 & 2.06\\
6 & 1 2 4 5 8 11   &  7 & 4 & 2.06\\
6 & 1 2 4 6 7 10   &  7 & 4 & 2.06\\
7 & 1 2 4 5 7 8 10   &  10 & 4 & 2.12\\
7 & 1 2 4 5 7 8 11   &  10 & 4 & 2.12\\
7 & 1 2 4 5 7 9 10   &  10 & 4 & 2.12\\
7 & 1 2 4 5 8 9 11   &  10 & 4 & 2.12\\
8 & 1 2 4 5 7 8 10 11   &  16 & 4 & 2.59\\
9 & 1 2 3 4 5 7 8 10 11   &  19 & 3 & 1.83\\
10 & 1 2 3 4 5 6 7 8 10 11   &  24 & 2 & 1.26\\
11 & 1 2 3 4 5 6 7 8 9 10 11   &  30 & 1 & 0.66\\
12 & 1 2 3 4 5 6 7 8 9 10 11 12   &  40 & 0 & 0\\
\hline
\end{tabular}
\label{table2}
\end{table}

\section{Musicality of some particular modes}

The following table lists some particular modes, along with their notes (all in the key of A), their musicality, and which of the other modes in the table they are equivalent.  The last column states the equivalences of the mode with the other modes' numbers, if any.  It lists "None" if it is not equivalent to any others in the table.

\bigskip

\begin{table}[h!]
\caption{Musicality of historic musical modes without the $6$ semitone interval}
\begin{tabular}{|l|l|l|l|l|}
\hline
Mode number &Mode name &Notes in the mode&Musicality&Equivalent  to\\
\hline
1 & Diminished & $A B C D D^{\#} F F^{\#} G^{\#} $ & 1.29 & 9  \\
2 & Enigmatic & $A A^{\#} C^{\#} D^{\#} F G G^{\#} $ & 0.85 &  None  \\
3 & Hungarian minor & $ A B C D^{\#} E F G^{\#}  $ & 1.11 & 11 13  \\
4 & Neapolitan & $A A^{\#} C D E F G^{\#} $  & 1.27 &  None  \\
5 & Pentatonic & $ A B C^{\#} E F^{\#}$  & 1.12 & 6 14 23  \\
6 & Pentatonic minor &$A C D E G$   & 1.12 & 5 14 23  \\
7 & Augmented & $A B C^{\#} D^{\#} F G$  & 0.88 &  None  \\
8 & Arabian & $A B C^{\#} D D^{\#} F G $  & 1.06 &  None  \\
9 & Diminished blues & $A A^{\#} C C^{\#} D^{\#} E F^{\#} G $   & 1.29 & 1  \\
10 & Balinese &$A A^{\#} C E F$  & 0.67 &  None  \\
11 & Byzantine &  $A A^{\#} C^{\#} D E F G^{\#} $  & 1.11 & 3 13  \\
12 & Chinese & $A C^{\#} D^{\#} E G^{\#}$   & 0.67 & 18 19  \\
13 & Double harmonic & $ A A^{\#} C^{\#} D E F G^{\#}$   & 1.11 & 3 11  \\
14 & Egyptian & A B D E G  & 1.12 & 5 6 23  \\
15& Eight tone Spanish & $A A^{\#} C C^{\#} D D^{\#} F G$& 1.29 &  None  \\
16 & Harmonic minor & $A B C D E F G^{\#}$& 1.27 & 27 28  \\
17 & Hindustan &$A B C^{\#} D E F G$  & 1.32 & 21 22  \\
18 & Hirajoshi & $ A C^{\#} D F^{\#} G^{\#} $   & 0.67 & 12 19  \\
19 & Japanese &$A A^{\#} D E F $ & 0.67 & 12 18  \\
20 & InSen & $A A^{\#} D E G $  & 0.45 &  None  \\
21 & Javanese & $A A^{\#} C D E F^{\#} G $   & 1.32 & 17 22  \\
22 & Melodic minor & $ A B C D E F^{\#} G^{\#}$  & 1.32 & 17 21  \\
23 & Mongolian & $A B C^{\#} E F^{\#} $  & 1.12 & 5 6 14  \\
24 & Pelog & $ A A^{\#} C D^{\#} G G^{\#} $ & 0.66 &  None  \\
25 & Persian &$A A^{\#} C^{\#} D D^{\#} F G^{\#}$  & 0.79 &  None  \\
26 & Prometheus &$A B C^{\#} D^{\#} F^{\#} G  $ & 1.1 &  None  \\
27 & Romanian minor & $A B C D{\#} E F^{\#} G$   & 1.27 & 16 28  \\
28 & Spanish Gypsy & $A A^{\#} C^{\#} D E F G $  & 1.27 & 16 27  \\

\hline
\end{tabular}
\label{table3} 
\end{table}
\newpage

The last table lists these modes again, with the alternate definition of a harmonic chord (from section 3).

\begin{table}[h!]
\caption{Musicality of historic musical modes with the $6$ semitone interval}
\begin{tabular}{|l|l|l|l|l|}
\hline
Mode number&Mode name&Notes in mode&Musicality&Equiv to\\
\hline

1 & Diminished & $A B C D D^{\#} F F^{\#} G^{\#} $ & 2.59 & 9  \\
2 & Enigmatic & $A A^{\#} C^{\#} D^{\#} F G G^{\#} $ & 1.06 &  None  \\
3 & Hungarian minor &$ A B C D^{\#} E F G^{\#}  $& 1.43 & 11 13  \\
4 & Neapolitan & $A A^{\#} C D E F G^{\#} $ & 1.48 &  None  \\
5 & Pentatonic &$ A B C^{\#} E F^{\#}$  & 1.12 & 6 14 23  \\
6 & Pentatonic minor & $A C D E G$  & 1.12 & 5 14 23  \\
7 & Augmented & $A B C^{\#} D^{\#} F G$  & 0.88 &  None  \\
8 & Arabian & $A B C^{\#} D D^{\#} F G $ & 1.32 &  None  \\
9 & Diminished Blues & $A A^{\#} C C^{\#} D^{\#} E F^{\#} G $ & 2.59 & 1  \\
10 & Balinese & $A A^{\#} C E F$  & 0.67 &  None  \\
11 & Byzantine & $A A^{\#} C^{\#} D E F G^{\#} $ & 1.43 & 3 13  \\
12 & Chinese & $A C^{\#} D^{\#} E G^{\#}$  & 0.67 & 18 19  \\
13 & Double harmonic &$ A A^{\#} C^{\#} D E F G^{\#}$  & 1.43 & 3 11  \\
14 & Egyptian & $A B D E G$  & 1.12 & 5 6 23  \\
15 & Eight tone Spanish & $A A^{\#} C C^{\#} D D^{\#} F G$  & 1.62 &  None  \\
16 & Harmonic minor & $A B C D E F G^{\#}$  & 2.12 & 27 28  \\
17 & Hindustan & $A B C^{\#} D E F G$  & 1.85 & 21 22  \\
18 & Hirajoshi &$ A C^{\#} D F^{\#} G^{\#} $ & 0.67 & 12 19  \\
19 & Japanese & $A A^{\#} D E F $ & 0.67 & 12 18  \\
20 & InSen & $A A^{\#} D E G $ & 0.89 &  None  \\
21 & Javanese & $A A^{\#} C D E F^{\#} G $ & 1.85 & 17 22  \\
22 & Melodic minor &$ A B C D E F^{\#} G^{\#}$ & 1.85 & 17 21  \\
23 & Mongolian & $A B C^{\#} E F^{\#} $ & 1.12 & 5 6 14  \\
24 & Pelog &$ A A^{\#} C D^{\#} G G^{\#} $ & 0.88 &  None  \\
25 & Persian & $A A^{\#} C^{\#} D D^{\#} F G^{\#}$  & 0.95 &  None  \\
26 & Prometheus & $A B C^{\#} D^{\#} F^{\#} G  $& 1.47 &  None  \\
27 & Romanian minor & $A B C D{\#} E F^{\#} G$  & 2.12 & 16 28  \\
28 & Spanish Gypsy & $A A^{\#} C^{\#} D E F G $ & 2.12 & 16 27  \\

\hline
\end{tabular}
\label{table4}
\end{table}

\end{document}